\documentclass[a4paper]{article}

\usepackage{amssymb, latexsym}

\usepackage{graphicx}

\usepackage{amsmath}
\usepackage{color}

\newtheorem{theorem}{Theorem}[section]

\newtheorem{corollary}[theorem]{Corollary}

\newtheorem{proposition}[theorem]{Proposition}

\newtheorem{remark}[theorem]{Remark}

\date{ }

\begin{document}
\title{\Large\bf $W$-entropy formulas on super Ricci flows 
and  Langevin deformation on Wasserstein space 
over Riemannian manifolds}
\author{Songzi Li\footnote{Supported by a Postdoctoral Fellowship at Beijing Normal University.}, \ \ \ Xiang-Dong Li\thanks{Research supported by NSFC No. 11371351 and Key Laboratory RCSDS, CAS, No. 2008DP173182.}\\
}

\maketitle

\thispagestyle{empty}
\begin{minipage}{120mm}
In this survey paper, we give an overview of our recent works on the study of the $W$-entropy for the heat equation associated with the Witten Laplacian on 
super-Ricci flows and the Langevin deformation on Wasserstein space over Riemannian manifolds. Inspired by Perelman's seminal work on the entropy formula for 
the Ricci flow, we prove 
the $W$-entropy formula for the heat equation associated with  the Witten Laplacian on $n$-dimensional complete Riemannian manifolds with the $CD(K, m)$-condition, and 
the $W$-entropy formula for the heat equation associated with the time dependent Witten Laplacian on $n$-dimensional  compact manifolds equipped with a $(K, m)$-super Ricci flow, where $K\in \mathbb{R}$ and $m\in [n, \infty]$.  Furthermore, we prove 
an analogue of the $W$-entropy formula for the geodesic flow on the Wasserstein space over Riemannian manifolds. Our result recaptures an important result due to Lott and Villani on 
the displacement convexity of  the Boltzmann-Shannon entropy on Riemannian manifolds with non-negative Ricci curvature.  To better understand the similarity between 
above two $W$-entropy formulas, we introduce the Langevin deformation of geometric flows on 
the cotangent bundle over the Wasserstein space and prove an extension of the $W$-entropy formula for the Langevin deformation. Finally, 
we make a discussion on  the $W$-entropy for the Ricci flow from the point of view of statistical mechanics and probability theory.  
\end{minipage}

\vskip1cm
\noindent{\it MSC2010 Classification}: primary 58J35, 58J65; secondary 60J60, 60H30.

\medskip

\noindent{\it Keywords}: $W$-entropy, Witten Laplacian,  Langevin deformation, $(K, m)$-super Ricci flows.


\section{Introduction}

Entropy was  introduced by R. Clausius  \cite{Cl1}  in 1865  in the study of thermodynamics.  In 1872, L. Boltzmann \cite{Btz1}  introduced the ${\rm H}$-entropy and formally derived the ${\rm H}$-theorem for the evolution equation of the probability distribution of ideal gas (now called the Boltzmann equation). The statistical interpretation of the {\rm H}-entropy was given by  Boltzmann \cite{Btz2}  in 1877. In 1948, C. Shannon \cite{Shan} introduced the Shannon entropy in the theory of communication and transformation of information.  In 1958, J. Nash \cite{Na} used 
the Boltzmann entropy to study the  continuity of solutions of parabolic and elliptic equations. 

Now entropy has been an important tool in many areas of mathematics. For example, the Kolmogorov-Sinai entropy plays an important role in the study of 
dynamical systems and ergodic theory, the exponential  decay of the Boltzmann entropy is closely related to the logarithmic Sobolev inequalities, the rate function in Sanov's theorem  in the theory of large deviation is the relative Boltzmann  entropy with respect to the reference measure.  More recently,  the displacement convexity of the Boltzmann-Shannon entropy 
or the R\'enyi  entropy is a key tool in Lott, Villani and Sturm's works \cite{LoV, Lo2, V1, V2, St1, St2, St3}  to develop analysis and geometry on metric measure spaces. 

In 1982, R. Hamilton \cite{H1} introduced the Ricci flow and initiated the program to prove the Poincar\'e conjecture and the Thurston geometrization conjecture using the Ricci flow (see also \cite{H2}). In a seminal paper \cite{P1}, Perelman gave a gradient flow reformulation for the Ricci flow and proved the $W$-entropy formula along the conjugate heat equation of the Ricci flow.  More precisely, 
let $M$ be an $n$-dimensional compact manifold and define
\begin{eqnarray*}
\mathcal{F}(g, f)=\int_M (R+|\nabla f|^2)e^{-f}dv,
\end{eqnarray*}
where $g\in \mathcal{M}=\{{\rm Riemannian\ metric \ on}\ M\}$,
$f\in C^\infty(M)$, $R$ denotes the scalar curvature on $(M, g)$,
and $dv$ denotes the volume measure on $(M, g)$.  Under the constraint condition
that the weighted volume measure 
\begin{eqnarray*}
d\mu=e^{-f}dv\end{eqnarray*}
is fixed, Perelman \cite{P1} proved that the gradient flow of $\mathcal{F}$ with respect to the standard $L^2$-metric on  $\mathcal{M}\times C^\infty(M)$ is given by the following modified Ricci flow for $g$ together with  the conjugate heat equation for $f$, i.e.,
\begin{eqnarray*}
\partial_t g&=&-2(Ric+\nabla^2 f),\\
\partial_t f&=&-\Delta f-R.
\end{eqnarray*}
Moreover, Perelman \cite{P1} introduced the $W$-entropy as follows
\begin{eqnarray}
W(g, f, \tau)=\int_M \left[\tau(R+|\nabla
f|^2)+f-n\right]{e^{-f}\over
 (4\pi\tau)^{n/2}}dv,\label{entropy-1}
 \end{eqnarray}
where $\tau>0$, and $f\in C^\infty(M)$ such that 
$$
\int_M (4\pi\tau)^{-n/2}e^{-f}dv=1,$$  and proved that if $(g(t), f(t), \tau(t))$ satisfies the evolution
equations
\begin{eqnarray}
\partial_t g&=&-2Ric,\label{RF}\\
\partial_t f&=&-\Delta f+|\nabla
f|^2-R+\frac{n} {2\tau},\label{r-c}\\
\partial_t \tau&=&-1,\nonumber
\end{eqnarray}
then the following remarkable 
$W$-entropy formula holds
\begin{eqnarray}
{d\over dt}W(g, f, \tau)=2 \int_M \tau\left|Ric+\nabla^2
f-{g\over 2\tau}\right|^2{e^{-f}\over (4\pi \tau)^{n/2}}dv.\label{Entropy-P}
\end{eqnarray}
In particular,  the $W$-entropy is monotonic  increasing in $t$
and the monotonicity is strict except that $(M, g(\tau), f(\tau))$
is a shrinking Ricci soliton, i.e.,
\begin{eqnarray*}
Ric+\nabla^2 f={g\over 2\tau}.\label{SRS}
\end{eqnarray*}
As an application,  Perelman \cite{P1} proved the no local
collapsing theorem, which ``removes the major stumbling block in
Hamilton's approach to geometrization'' and plays an important role in the final resolution of the
Poincar\'e conjecture and Thurston's geometrization conjecture. 

It is natural and interesting to ask the problems what is the hidden idea for Perelman to introduce the mysterious
$W$-entropy,  what is the reason for him to call the quantity in $(\ref{entropy-1})$ the $W$-entropy, and whether there is 
some essential link between the $W$-entropy 
and the Boltzmann entropy in statistical mechanics and probability theory. 

Inspired by Perelman \cite{P1} and related works \cite{N1, N2},  the second author of this paper proved in \cite{Li07}  the $W$-entropy formula for the heat equation of the Witten Laplacian on compact
 Riemannian manifolds with the $CD(0, m)$-condition and gave a probabilistic interpretation of the $W$-entropy for the Ricci flow. Later,   the $W$-entropy formula and a rigidity theorem for the $W$-entropy were proved in \cite{Li12, Li13} for the fundamental solution to the heat equation of the Witten Laplacian on complete 
 Riemannian manifolds with the $CD(0, m)$-condition,  and  the $W$-entropy formula was proved in \cite{Li11} for the Fokker-Planck equation of the 
 Witten Laplacian on complete 
 Riemannian manifolds with the $CD(0, m)$-condition.  The relationship between Perelman's $W$-entropy formula for the Ricci flow and the Boltzmann 
 ${\rm H}$-theorem for the Boltzmann equation was discussed in \cite{Li12b}  from the point of view of the statistical mechanics.  In  \cite{LL15a, LL15b,  LL17a, LL17b, SLi15},  we extended the $W$-entropy formula to the heat equation of the Witten Laplacian on complete Riemannian manifolds with the $CD(K, m)$-condition and on compact Riemannian manifolds equipped 
with $(K, m)$-super Ricci flows. Moreover, we proved in \cite{LL16, SLi15} an analogue of the $W$-entropy formula for the geodesic flow 
on the Wasserstein space over Riemannian manifolds with the $CD(0, m)$-condition, 
which recaptures an important result due to Lott and Villani \cite{LoV, Lo2} on the displacement convexity of the Boltzmann-Shannon entropy on the 
Wasserstein space over Riemannian manifolds with non-negative Ricci curvature. To better understand the similarity between the $W$-entropy formula for the heat equation of the Witten Laplacian and the $W$-entropy formula for the geodesic flow on the Wasserrstein space over Riemannian manifolds, we introduced in \cite{LL16, SLi15} the Langevin deformation of geometric flows on the Wasserstein space over Riemannian manifolds, which interpolates the backward gradient flow of the Boltzmann-Shannon entropy and the geodesic flow on the Wasserstein space, and proved an extension of the $W$-entropy formula for the Langevin deformation. The rigidity models are also proposed for the Langevin deformation of flows.  
In particular, two rigidity theorems were proved in \cite{Li12, Li13, LL16, SLi15}  for the gradient flow of the Boltzmann-Shannon entropy and the geodesic flow on the Wasserstein space over complete Riemannian manifolds with the  $CD(0, m)$-condition.

The purpose of this survey paper is to give an overview of  our  works in \cite{Li07, Li12, Li12b, Li13, LL15a, LL15b, LL16,  LL17a, LL17b, SLi15} and to make a 
discussion on the $W$-entropy for the Ricci flow from the point of view of  statistical mechanics and probability theory. 

In 2016, the second author of this paper was invited to give a  Special Invited Talk in the 2016 
Autumn Meeting of Mathematical Society of Japan. This paper is an improved version of the abstract for this meeting. He would like to thank the committee members of the Mathematical Society of Japan, in particular, Professor S. Aida and Professor K. Kuwae,  for their interests and invitation. We would like to thank Professor Feng-Yu Wang 
for inviting us to submit this paper as a Special Invited Paper to SCIENCE CHINA Mathematics and the Mathematical Society of Japan for their permission.

\section{$W$-entropy formulas for Witten Laplacian on Riemannian manifolds}

Since Perelman's preprint \cite{P1} was published on Arxiv in 2002, many
people have studied the $W$-like entropy for other geometric
flows on Riemannian manifolds \cite{N1, N2, Ec, LNVV, KN}. Let
$(M, g)$ be an $n$-dimensional complete Riemannian manifold with a fixed metric (and with bounded geometry condition), and let
$$u={e^{-f}\over (4\pi t)^{n/2}}$$
be a positive solution to the linear heat equation
\begin{eqnarray}
\partial_t u=\Delta u \label{heat-1}
\end{eqnarray}
with $\int_M u(x, 0)dv(x)=1$. In \cite{N1, N2}, Ni introduced the $W$-entropy  for the
linear heat equation $(\ref{heat-1})$ by
\begin{eqnarray}
W(f, t)=\int_M \left[t |\nabla
f|^2+f-n\right]{e^{-f}\over
 (4\pi t)^{n/2}}dv,\label{entropy-2}
\end{eqnarray}
and proved the following $W$-entropy formula
\begin{eqnarray}
{d\over dt}W(f, t)&=&-2\int_M t\left(
\left|\nabla^2 f-{g\over 2t}\right|^2+Ric(\nabla f, \nabla
f)\right){e^{-f}\over (4\pi t)^{n/2}}dv. \label{NiW}
\end{eqnarray}
In particular,  the $W$-entropy for the
linear heat equation $(\ref{heat-1})$ is decreasing on complete Riemannian
manifolds with non-negative Ricci curvature. In \cite{LX}, Li and Xu extended Ni's $W$-entropy formula $(\ref{NiW})$ to the heat equation $\partial_t u=\Delta u$ on complete Riemannian manifolds with fixed metric satisfying  $Ric\geq -Kg$, where $K\geq 0$ is a constant.

Let $(M, g)$ be a complete Riemannian manifold, $\phi\in C^2(M)$. Let  $d\mu=e^{-\phi}dv$, where $dv$ is the Riemannian volume measure on $(M, g)$. The Witten Laplacian, called also the weighted Laplacian,
\begin{eqnarray*}
L =\Delta -\nabla \phi\cdot\nabla \label{WL}
\end{eqnarray*}
is a self-adjoint and non-negative operator on $L^2(M, \mu)$. By It\^o's calculus, one can construct the symmetric diffusion process $X_t$ associated to the Witten Laplacian by solving the SDE
\begin{eqnarray*}
dX_t=\sqrt{2}dW_t-\nabla\phi(X_t)dt,
\end{eqnarray*}
where $W_t$ is the Brownian motion on $M$. Moreover, it is well known that the transition probability density of the diffusion process $X_t$ is  the fundamental solution to the heat equation of $L$, i.e., the heat kernel of the Witten Laplacian $L$. In view of this, it is a fundamental problem to study the heat equation and related properties for the Witten Laplacian on Riemannian 
manifolds with various geometric conditions.

To develop the study of the $W$-entropy formula for the heat equation of the Witten Laplacian, we need to introduce some notations. Let $n={\rm dim}~ M$, and $m\geq n$ 
a constant.  Following Bakry and Emery \cite{BE}, we introduce 
\begin{eqnarray*}
Ric_{m, n}(L)=Ric+\nabla^2\phi-{\nabla\phi\otimes\nabla\phi\over m-n}, 
\end{eqnarray*}
and call it  the $m$-dimensional Bakry-Emery Ricci curvature associated with the Witten Laplacian $L$ on $(M, g, \phi)$.  
We make the convention that $m=n$ if and only if $L=\Delta$ and $\phi$ is a constant. In this case,  $Ric_{m, n}(L)=Ric$. When $m=\infty$,  we introduce 
\begin{eqnarray*}
Ric(L)=Ric+\nabla^2\phi.
\end{eqnarray*}
Following Bakry and Emery \cite{BE}, we say that the $CD(K, m)$ condition holds for the Witten Laplacian $L=\Delta-\nabla\phi\cdot\nabla$ on $(M, g, \phi)$ 
if and only if $Ric_{m, n}(L)\geq K$, where $K\in \mathbb{R}$ and $m\in [n, \infty]$.

\subsection{The case of $CD(0, m)$-condition} 

In \cite{Li07, Li12, Li13},  inspired by Perelman's work on the $W$-entropy formula for Ricci flow, the second author of this paper proved the $W$-entropy formula for the heat equation associated with the Witten-Laplacian on complete Riemannian manifolds with the $CD(0, m)$-condition, which extends the above mentioned result due to Ni \cite{N1, N2}.  More precisely, we have

\begin{theorem}\label{Th-A} (\cite{Li07, Li12, Li13})  Let $(M, g)$ be a compact Riemannian manifold, or a complete Riemannian manifold with bounded geometry condition\footnote{Here we say that $(M, g)$ satisfies the bounded geometry condition if the Riemannian curvature tensor ${\rm Riem}$ and its covariant derivatives $\nabla^k {\rm Riem}$ are uniformly bounded on $M$, $k=1, 2, 3$.}, and $\phi\in C^4(M)$ with $\nabla\phi\in C_b^3(M)$. Let $m\in [n, \infty)$, and  $u={e^{-f}\over (4\pi t)^{m/2}}$ be a positive solution of the heat equation $\partial_t u=Lu$ when $M$ is compact, or the fundamental solution associated with the Witten Laplacian, i.e., the heat kernel to
the heat equation $\partial_t u=Lu$, when $M$ is complete non-compact.
Let 
\begin{eqnarray*}
H_m(u, t)=-\int_M u\log u d\mu-{m\over 2}(1+\log(4\pi t)).
\end{eqnarray*}
Define the $W$-entropy for the Witten-Laplacian by
\begin{eqnarray*}
W_m(u, t)={d\over dt}(tH_m(u)).
\end{eqnarray*}
Then
\begin{eqnarray}
W_m(u, t)=\int_M\left[t|\nabla f|^2+f-m\right]{e^{-f}\over (4\pi t)^{m/2}}d\mu,\label{WW-0}
\end{eqnarray}
and
\begin{eqnarray}
{d\over dt} W_m(u, t)&=&-2\int_M
t \left(\left|\nabla^2 f-{g\over 2t}\right|^2+Ric_{m,
n}(L)(\nabla f,
\nabla f)\right)ud\mu\nonumber\\
& & \hskip2cm  -{2\over m-n}\int_M t \left({\nabla\phi\cdot\nabla
f}+{m-n\over 2t}\right)^2ud\mu.\label{W-1}
\end{eqnarray}
\end{theorem}

By calculation and integration by part, we have 
\begin{eqnarray*}
{d\over dt}H_m(u, t)=-\int_M \left(L\log u+{m\over 2t}\right)ud\mu.\label{Hentropy2}
\end{eqnarray*} 
By \cite{Li05, Li12},  if $Ric_{m, n}(L)\geq 0$, the generalized Li-Yau Harnack inequality (\cite{LY}) holds 
\begin{eqnarray*}
L\log u+{m\over 2t}\geq 0, \ \ \ \forall t>0.\label{LY-1}
\end{eqnarray*}
Therefore, $H_m(u, t)$ is non-increasing in time $t$ for the heat equation $\partial_t u=Lu$ on complete Riemannian manifolds with the $CD(0, m)$-condition, 
i.e., $Ric_{m, n}(L)\geq 0$. 

As a corollary of Theorem \ref{Th-A},  if $(M, g, \phi)$ is complete Riemannian manifold with the bounded geometry condition and $Ric_{m, n}(L)\geq 0$, then the $W$-entropy for the heat equation $\partial_t u=Lu$ is decreasing in time $t$, i.e.,
\begin{eqnarray*}
{d\over dt}W_m(u, t)\leq 0, \ \ \ \ \ \ \ \forall t\geq 0.
\end{eqnarray*}
Moreover, under the condition $Ric_{m, n}(L)\geq 0$,  it was proved in \cite{Li12}  that $W_m(u, t)$ attains its minimum at some point $t=t_0>0$, i.e.,
\begin{eqnarray*}
{d\over dt}W_m(u, t)=0 \ \ \ \ {\rm at\ some}\ \ t=t_0>0,
\end{eqnarray*}
if and only if $(M, g)$ is isometric to Euclidean space $\mathbb{R}^n$, $m=n$, $\phi\equiv C$ for a constant $C\in \mathbb{R}$, and
\begin{eqnarray*}
u(x, t)={e^{-{\|x\|^2\over 4t}}\over (4\pi t)^{n/2}}, \ \ \ \ \forall x\in \mathbb{R}^n, t>0.
\end{eqnarray*}
In other words, the Euclidean space $\mathbb{R}^n$ is the unique equilibrium state for the $W$-entropy of the Witten-Laplacian in the statistical ensemble of 
complete Riemannian manifolds $(M, g, \phi)$  with bounded geometry condition and with the $CD(0, m)$-condition.

In \cite{LL15a},  we gave a new proof of Theorem \ref{Th-A} by using the warped product approach.  Let $m\in \mathbb{N}$, $m\geq n$. Let
$\widetilde{M}=M\times N$, where $(N, g_N)$ is a compact Riemannian manifold with dimension $m-n$. Let $\phi\in C^2(M)$. We consider the following warped product metric on $\widetilde{M}$:
\begin{eqnarray*}
\widetilde{g}=g_M\bigoplus e^{-{2\phi\over m-n}}g_N.\label{WPM}
\end{eqnarray*}
Let  $\nu_N$ be the normalized volume measure on $N$. Then the volume measure on $(\widetilde{M}, \widetilde{g})$ is
$$dvol_{\widetilde{M}}=d\mu\otimes d\nu_N.$$
Let $\widetilde\nabla$ be the Levi-Civita connection on $(\widetilde{M}, \widetilde{g})$,  $\widetilde{\nabla}^2$ and $\widetilde{\Delta}$ 
the Hessian and the Laplace-Beltrami operator on  $(\widetilde{M}, \widetilde{g})$.  By direct calculation, we have 
\begin{eqnarray}
\left|\widetilde{\nabla}^2 f-{\widetilde{g}\over 2t}\right|^2
=\left|\nabla^2 f-{g\over 2t}\right|^2+{1\over m-n}\left({\nabla\phi\cdot\nabla f}+{m-n\over 2t}\right)^2,\label{ccc}
\end{eqnarray}
and
\begin{eqnarray*}
\widetilde{\Delta}=L+e^{-{2\phi\over m-n}}\Delta_{N}.
\end{eqnarray*}

\medskip
\noindent
{\bf A new proof of Theorem \ref{Th-A} (\cite{LL15a})}.\ \ To avoid technical issue, we only consider the case of compact manifolds. Let 
$u={e^{-f}\over (4\pi t)^{m/2}}$ be a positive solution to the heat equation $\partial_t u=Lu$. Then it satisfies the following heat equation on $(\widetilde{M}, \widetilde{g})$
\begin{eqnarray*}
\partial_t u=\widetilde{\Delta} u.
\end{eqnarray*}
Since $f$ depends only on the variable in the $M$-direction, we have $\widetilde\nabla f=\nabla f$. Therefore the $W$-entropy $W_m(u, t)$ defined by 
$(\ref{WW-0})$ coincides with the $W$-entropy $\widetilde{W}_{m}(u, t)$ defined on $(\widetilde{M}, \widetilde{g})$ as follows
\begin{eqnarray*}
\widetilde{W}_m(u, t)=\int_{\widetilde{M}}\left[t|\widetilde\nabla f|^2+f-m\right]{e^{-f}\over (4\pi t)^{m/2}}dvol_{\widetilde{M}}.
\end{eqnarray*}
Applying Ni's $W$-entropy formula  $(\ref{NiW})$ to $(\widetilde{M}, \widetilde{g})$,  we have
\begin{eqnarray}
{d\over dt}\widetilde{W}_m(u, t)&=&-2\int_{\widetilde{M}} t \left(\left|\widetilde\nabla^2 f-{\widetilde g\over 2t}\right|^2
+\widetilde{Ric}(\widetilde\nabla f, \widetilde\nabla f)\right)ud\mu dv_N.\label{WP-W-1}
\end{eqnarray}
By $(\ref{ccc})$ and $\widetilde{Ric}(\widetilde\nabla f, \widetilde\nabla f)=Ric_{m, n}(L)(\nabla f,
\nabla f)$, we derive $(\ref{W-1})$ from $(\ref{WP-W-1})$.  \hfill $\square$

\begin{remark}\label{rem1} {\rm One of the advantages of the above proof is that: when $m\in \mathbb{N}$ and $m>n$, the quantity ${1\over {m-n}}\left(\nabla \phi\cdot \nabla f+{m-n\over 2t}\right)^2$ appeared in the $W$-entropy formula in Theorem \ref{Th-A} has a natural geometric interpretation. It corresponds to the vertical component of the  quantity $\left|\widetilde\nabla^2 f-{\widetilde{g}\over 2t}\right|^2$  on $\widetilde{M}=M\times N$ equipped with the warped product metric $(\ref{WPM})$. 
}
\end{remark}

\subsection{The case of $CD(K, m)$-condition} 

Theorem \ref{Th-A} can be viewed as the $W$-entropy formula for the heat equation of the Witten Laplacian on complete Riemannian manicolds with the $CD(0, m)$-condition.  It is natural to raise the question whether we can extend Theorem \ref{Th-A}  to 
the heat equation of the Witten Laplacian  on complete Riemannian manifolds with the $CD(K, m)$-condition for general $K\in \mathbb{R}$ and $m\in [n, \infty]$. 

In  \cite{LL15a},  we extended Theorem \ref{Th-A}  to the  Witten Laplacian
on complete Riemannian manifolds with the $CD(K, m)$-condition for $K\in \mathbb{R}$ and $m\in [n, \infty)$.

\begin{proposition} (\cite{LL15a})\label{WW1} Let $(M, g)$ be a complete Riemannian manifold with bounded geometry condition, and $\phi\in C^4(M)$ satisfying the condition in Theorem \ref{Th-A}.
 Let $u$ be a positive solution to the heat equation $\partial_t u=Lu$. Then, under  the $CD(-K, m)$-condition, i.e., $Ric_{m, n}(L)\geq -K$, where $K\in \mathbb{R}$ and $m\in [n, \infty)$,  the following Harnack inequality holds 
\begin{eqnarray*}
{|\nabla u|^2\over u^2}-\left(1+{2\over 3}Kt\right){\partial_t u\over u}\leq {m\over 2t}+{mK\over 2}\left(1+{Kt\over 3}\right), \ \ \ \forall t>0. \label{HIKm}
\end{eqnarray*}
\end{proposition}

\begin{theorem} (\cite{LL15a, LL15b}) \label{Th-B} Let $u=\frac{e^{-f}}{(4\pi t)^{m/2}}$ be the fundamental solution to the heat equation $\partial_t u=Lu$. Under the same assumption as in Theorem \ref{Th-A}, define 
\begin{eqnarray}
H_{m,K}(u, t)=-\int_M u\log u d\mu-{m\over 2}(1+\log(4\pi t))-\frac m2Kt\Big(1+\frac16Kt\Big), \label{HKm}
\end{eqnarray}
and the $W$-entropy by the Boltzmann formula
\begin{eqnarray}
W_{m, K}(u, t)={d\over dt}(tH_{m,K}(u)). \label{WKm}
\end{eqnarray}
Then
\begin{eqnarray*}
W_{m, K}(u, t)=\int_M\left(t|\nabla f|^2+f-m\Big(1+\frac12Kt\Big)^2\right)u d\mu,\label{WmK}
\end{eqnarray*}
and
\begin{eqnarray}
& &\frac{d}{dt}W_{m, K}(u, t)+2t\int_M\left(\Big|\nabla^2 f-\left(\frac1{2t}+\frac K2\right) g\Big|^2\right)u d\mu\nonumber\\
& &\hskip2cm +{2t\over m-n}\int_M \left(\nabla \phi\cdot\nabla f+(m-n)\Big(\frac1{2t}+\frac K2\Big)\right)^2u\ d\mu\nonumber\\
& &\hskip3cm  =2t\int_M ({\rm Ric}_{m,n}(L)+Kg)(\nabla f, \nabla f)u\ d\mu.\label{WmK20}
\end{eqnarray}
In particular, if the $CD(-K, m)$-condition holds, i.e., $Ric_{m, n}(L)\geq -K$, then
\begin{eqnarray*}
{d\over dt}H_{m, K}(u, t)\leq 0,
\end{eqnarray*}
and
\begin{eqnarray*}
\frac{d}{dt}W_{m, K}(u, t)
\leq 0.
\end{eqnarray*}
Moreover, under the $CD(-K, m)$-condition, the left hand side  in $(\ref{WmK20})$  equals to zero at some $t=t_0>0$  if and only if $(M, g, \phi)$ is a $(-K, m)$-Ricci soliton, i.e., 
\begin{eqnarray*}
Ric_{m, n}(L)=-Kg.
\end{eqnarray*}
\end{theorem}

\subsection{The case of $CD(K, \infty)$-condition}

When $m=\infty$,  we cannot use the above definition formulas $(\ref{HKm})$  and $(\ref{WKm})$ to introduce $H_{K, \infty}(u, t)$ and to define the $W$-entropy for the Witten Laplacian on Riemannian manifolds with the $CD(K, \infty)$-condition. 
Based on the reversal logarithmic Sobolev inequality due to Bakry and Ledoux \cite{BL}, we proved the following result. 

\begin{theorem} (\cite{LL15b,  LL17b, SLi15}) \label{Th-C} Let $M$ be a complete Riemannian manifold with
bounded geometry condition, $\phi\in C^4(M)$ with $\nabla\phi\in
C_b^3(M)$. Suppose that $Ric+\nabla^2\phi\geq K$, where $K\in \mathbb{R}$ is
a constant. Let $u(\cdot, t)=P_tf$ be a positive solution to the heat
equation $\partial_t u=Lu$ with $u(\cdot, 0)=f$, where $f$ is a positive and measurable function on $M$. Let
\begin{eqnarray*}
H_{K}(f, t)=D_K(t)\int_M (f\log f-P_tf\log P_tf )d\mu,
\end{eqnarray*}
where $D_0(t)={1\over t}$ and $D_{K}(t)={2K\over 1-e^{-2Kt}}$ for $K\neq 0$.Then for all  $t>0$
\begin{eqnarray*}
{d\over dt}H_{K}(f, t)\leq 0,
\end{eqnarray*}
and for all  $t>0$, we have
\begin{eqnarray*}
{d^2\over dt^2}H_K(t)+2K\coth(Kt) {d\over dt}H_K(t)
\leq - 2D_K(t)\int_M |\nabla^2\log P_tf|^2P_tfd\mu. \label{KK2}
\end{eqnarray*}
\end{theorem}

Theorem \ref{Th-C}  suggests us a new way to introduce the $W$-entropy  for the Witten Laplacian on Riemannian manifolds with the $CD(K, \infty)$-condition
\begin{eqnarray*}
W_K(f, t)=H_K(f, t)+{\sinh(2Kt)\over 2K}{d\over dt}H_K(f, t).
\end{eqnarray*}
In this way,  for all $t>0$, we prove that (\cite{LL15b})
\begin{eqnarray}
& &{d\over dt}W_{K}(f, t)+(e^{2Kt}+1)\int_M  |\nabla^2
\log P_tf|^2 P_tf d\mu\nonumber\\
& & \hskip1cm =-(e^{2Kt}+1)\int_M (Ric(L)-Kg)(\nabla\log P_tf, \nabla\log P_tf)P_tfd\mu.\label{WW1}
\end{eqnarray}
In particular, for all $t>0$, we have
\begin{eqnarray*}
{d\over dt}W_{K}(f, t)\leq 0.
\end{eqnarray*}
Moreover, under the $CD(K, \infty)$-condition, the left hand side in  $(\ref{WW1})$ equals to zero at some $t=t_0>0$ if and only if $(M, g, \phi)$ is a $K$-Ricci soliton, i.e., 
\begin{eqnarray*}
Ric+\nabla^2\phi=Kg.
\end{eqnarray*}

\section{$W$-entropy formulas for Witten Laplacian on $(K, m)$-super Ricci flows}

In Section $2$, we extend the $W$-entropy formula to the heat equation of the Witten Laplacian on complete Riemannian manifolds with the $CD(K, m)$-condition. It is 
an interesting question whether we can further extend the $W$-entropy formula to the heat equation associated with  the 
time dependent Witten Laplacian on compact or complete Riemannian manifolds with time dependent metrics and potentials.

We now introduce the notion of the $(K, m)$-super Ricci flow on manifolds with time dependent metrics and potentials. By definition, 
we call $(M, g(t), \phi(t), t\in [0, T])$ a $(K, m)$-super Ricci flow if 
\begin{eqnarray*}
{1\over 2}{\partial g\over \partial t}+Ric_{m, n}(L)\geq Kg, \ \ \ \ \forall ~ t\in [0, T],
\end{eqnarray*}
where $K\in \mathbb{R}$ and $m\in [n, \infty]$ are two constants. 
Note that  a Riemannian manifold equipped with a stationary $(K, m)$-super Ricci flow (i.e., $(g(t), \phi(t))$ is independent of time) if and only if the $CD(K, m)$-condition holds, i.e., 
\begin{eqnarray*}
Ric_{m, n}(L)\geq Kg.
\end{eqnarray*} 
In the case $m=n$, the notion of the $(K, n)$-super Ricci flow is indeed the $K$-super Ricci flow in geometric analysis
\begin{eqnarray*}
{1\over 2}{\partial g\over \partial t}+Ric\geq Kg, \ \ \ \ \forall ~ t\in [0, T], 
\end{eqnarray*}
and in the case $m=\infty$, the $(K, \infty)$-super Ricci flow equation reads
\begin{eqnarray*}
{1\over 2}{\partial  g\over \partial t}+Ric(L)\geq Kg, \ \ \ \ \forall ~ t\in [0, T]. \label{KKK}
\end{eqnarray*}
In view of this,  the Perelman Ricci flow is indeed the $(0, \infty)$-Ricci flow together with the conjugate heat equation
\begin{eqnarray*}
{\partial g\over \partial t}&=&-2Ric(L),\\
 \ \ \ {\partial\phi\over \partial t}&=&{1\over 2}{\rm Tr} \left({\partial g\over \partial t}\right).
\end{eqnarray*}

Let $(M, g(t), \phi(t), t\in [0, T])$ be a complete Riemannian manifold with a family of time dependent metrics $g(t)$ and potentials $\phi(t)$. Let
$$L=\Delta_{g(t)}-\nabla_{g(t)}\phi(t)\cdot\nabla_{g(t)}$$
be the time dependent Witten Laplacian on $(M, g(t), \phi(t))$. Let
$$d\mu(t)=e^{-\phi(t)}dvol_{g(t)}.$$
Suppose that 
\begin{eqnarray}
{\partial \phi\over \partial t}={1\over 2}{\rm Tr}\left( {\partial g\over
\partial t}\right).\label{conjugate}
\end{eqnarray}
Then $\mu(t)$ is independent of $t\in [0, T]$, i.e.,
\begin{eqnarray*}
{\partial d\mu(t)\over \partial t}=0, \ \ \ t\in [0, T].
\end{eqnarray*}

We now state the main results of this section, which extend
Theorems~\ref{Th-A}, \ref{Th-B} and  \ref{Th-C} to the heat equation associated with the time dependent
Witten Laplacian on compact manifolds with a $(K, m)$-super Ricci flow, where $K\in \mathbb{R}$ and $m\in [n, \infty]$.

\subsection{The case of $(0, m)$-super Ricci flow}

In \cite{LL15a}, we proved the $W$-entropy formula to the heat equation associated with  the time dependent Witten Laplacian on compact manifolds equipped with a $(0, m)$-super Ricci flow, which 
can be regarded as the $m$-dimensional analogue of Perelman's $W$-entropy formula for the Ricci flow.

\begin{theorem}\label{Th-D} (\cite{LL15a})  Let $(M, g(t), \phi(t), t\in [0, T])$ be a compact manifold with family of time dependent metrics and $C^2$-potentials. Suppose that $g(t)$ and $\phi(t)$
satisfy the conjugate equation $(\ref{conjugate})$.
Let  $u={e^{-f}\over (4\pi t)^{m/2}}$ be a positive solution of the heat equation
\begin{eqnarray*}
\partial_t u = Lu
\end{eqnarray*}
with initial data $u(0)$ satisfying $\int_M
u(0)d\mu(0)=1$.
Let
\begin{eqnarray*}
H_m(u, t)=-\int_M u\log u d\mu-{m\over 2}(1+\log(4\pi t)).
\end{eqnarray*}
Define
\begin{eqnarray*}
W_m(u, t)={d\over dt}(tH_m(u)).
\end{eqnarray*}
Then
\begin{eqnarray*}
W_m(u, t)=\int_M \left[t|\nabla f|^2+f-m\right]ud\mu,
\end{eqnarray*}
and
\begin{eqnarray}
& &{d\over dt}W_m(u, t)+2t\int_M \left|\nabla^2 f-{g\over 2t}\right|^2ud\mu+{2t\over m-n}\int_M \left(\nabla \phi\cdot \nabla f+{m-n\over 2t}\right)^2  ud\mu\nonumber\\
& &\hskip3cm =-2t\int_M \left({1\over 2}{\partial g\over \partial t}+Ric_{m, n}(L)\right)(\nabla f, \nabla f)ud\mu.\label{NW}
\end{eqnarray}
In particular, if $\{g(t), \phi(t), t\in (0, T]\}$ is a $(0, m)$-super Ricci flow and satisfies the conjugate equation $(\ref{conjugate})$, then $W_m(u, t)$ is decreasing in $t\in (0, T]$, i.e.,
\begin{eqnarray*}
{d\over dt}W_m(u, t)\leq 0, \ \ \ \forall t\in (0, T].
\end{eqnarray*}
Moreover, the left hand side in $(\ref{NW})$ identically  equals to zero  on $(0, T]$ if and only if $(M, g(t), \phi(t), t\in (0, T])$ is a $(0, m)$-Ricci flow in the sense that
\begin{eqnarray*}
{\partial g\over \partial t}&=&-2{\rm Ric}_{m,n}(L),\\
{\partial \phi\over \partial t}&=&{1\over 2} {\rm Tr}\left( {\partial g\over \partial t}\right).
\end{eqnarray*}

\end{theorem}

\subsection{The case of $(K, m)$-super Ricci flow}

In general we have the following result which extends Theorem \ref{Th-B} to $(K, m)$-super Ricci flow for general $K\in \mathbb{R}$ and $m\in [n, \infty)$.

\begin{theorem}\label{Th-E} (\cite{LL15a, LL15b}) Under the same notation as in Theorem \ref{Th-D},  define 
\begin{eqnarray}
H_{m,K}(u, t)=-\int_M u\log u d\mu-{m\over 2}(1+\log(4\pi t))-\frac m2Kt\Big(1+\frac16Kt\Big),  \label{HmK}
\end{eqnarray}
and
\begin{align}
W_{m,K}(u, t)={d\over dt}(tH_{m,K}(u)). \label{WmK}
\end{align}
Then
\begin{eqnarray*}
W_{m,K}(u, t)=\int_M\left[t|\nabla f |^2+f-m\Big(1+\frac12Kt\Big)^2\right]ud\mu,\label{WmK-0}
\end{eqnarray*}
and
\begin{align}\label{WMK}
& {d\over dt}W_{m,K}(u, t)+2 t\int_M\Big|\nabla^2 f-\left(\frac1{2t}+\frac{K}{2}\right)g\Big|^2u d\mu\nonumber\\
&\ \ \ \ +\frac{2t}{m-n}\int_M\left(\nabla \phi\cdot\nabla f+(m-n)\Big(\frac1{2t}+\frac K2\Big)\right)^2u d\mu\nonumber\\
&\ \ \ \ \ =-2 t\int_M\left({1\over 2}{\partial g\over \partial t}+{\rm Ric}_{m,n}(L)+Kg\right)(\nabla f, \nabla f) ud\mu.
\end{align}
In particular, if $(M, g(t), \phi(t), t\in (0, T])$ is a $(-K, m)$-super Ricci flow and satisfies the conjugate equation $(\ref{conjugate})$, then $W_{m,K}(u, t)$ is decreasing in $t\in (0, T]$, i.e.,
\begin{eqnarray*}
{d\over dt}W_{m,K}(u, t)\leq 0, \ \ \ \forall t\in (0, T].
\end{eqnarray*}
Moreover, the left hand side in $(\ref{WMK})$ identically  equals to zero on $(0, T]$ if and only if $(M, g(t), \phi(t), t\in (0, T])$ is a $(-K, m)$-Ricci flow in the sense that
\begin{eqnarray*}
{\partial g\over \partial t}&=&-2({\rm Ric}_{m,n}(L)+Kg),\\
{\partial \phi\over \partial t}&=&{1\over 2} {\rm Tr}\left( {\partial g\over \partial t}\right).
\end{eqnarray*}
\end{theorem}

\subsection{The case of $(K, \infty)$-super Ricci flow}

In \cite{LL15b, LL17b, SLi15}, we proved the equivalence between the $(K, \infty)$-super Ricci flow and two families of 
 logarithmic Sobolev inequalities for the time dependent Witten Laplacian on Riemannian manifolds with time dependent metrics and potentials. Based on this result, 
we have the following $W$-entropy formula  for the time dependent
Witten Laplacian on compact Riemannian manifolds with $(K, \infty)$-super Ricci flow, which can be viewed as the natural extension of the $W$-entropy formula for the heat equation of the Witten Laplacian on 
complete Riemannian manifolds with the $CD(K, \infty)$-condition.

\begin{theorem}\label{Th-F} (\cite{LL15b, LL17b, SLi15}) Let $(M, g(t), \phi(t), t\in [0, T])$ be a compact $(K, \infty)$-super Ricci flow satisfying the conjugate heat equation
$(\ref{conjugate})$. Let $u(\cdot, t)=P_tf$ be a positive solution to the heat
equation $\partial_t u=Lu$ with $u(\cdot, 0)=f$, where $f$ is a positive and measurable function on $M$. Define
\begin{eqnarray*}
H_{K}(f, t)=D_K(t)\int_M (f\log f-P_tf\log P_tf )d\mu,
\end{eqnarray*}
where $D_0(t)={1\over t}$ and $D_{K}(t)={2K\over 1-e^{-2Kt}}$ for $K\neq 0$.Then for all $t\in [0, T]$
\begin{eqnarray*}
{d\over dt}H_{K}(f, t)\leq 0,
\end{eqnarray*}
and for all $t\in (0, T]$, we have
\begin{eqnarray*}
{d^2\over dt^2}H_K(t)+2K\coth(Kt) {d\over dt}H_K(t)
\leq - 2D_K(t)\int_M |\nabla^2\log P_tf|^2P_tfd\mu.
\end{eqnarray*}
Define the $W$-entropy by the revised Boltzmann entropy formula
\begin{eqnarray*}
W_K(f, t)=H_K(f, t)+{\sinh(2Kt)\over 2K}{d\over dt}H_K(f, t).
\end{eqnarray*}
Then for all  $t\in (0, T]$, we have
\begin{eqnarray}
& &{d\over dt}W_{K}(f, t)+(e^{2Kt}+1)\int_M  |\nabla^2
\log P_tf|^2 P_tf d\mu\nonumber\\
& &\hskip0.5cm =-(e^{2Kt}+1)\int_M \left({1\over 2}{\partial g\over \partial t}+Ric(L)-Kg\right)(\nabla\log P_tf, \nabla\log P_tf)P_tfd\mu.\label{WW2}
\end{eqnarray}
In particular, for all $t\in (0, T]$, we have
\begin{eqnarray*}
{d\over dt}W_{K}(f, t)\leq 0.
\end{eqnarray*}
Moreover, the left hand side of $(\ref{WW2})$  identically equals to zero on $(0, T]$ if and only if $(M, g(t), \phi(t))$ is the $(K, \infty)$-Ricci flow 
satisfying the conjugate equation $(\ref{conjugate})$, i.e., for all $t\in (0, T]$, 
\begin{eqnarray*}
{\partial g\over \partial t}&=&-2(Ric+\nabla^2\phi-Kg), \\
 {\partial\phi\over \partial t}&=&-R-\Delta \phi+nK.
\end{eqnarray*}

\end{theorem}

\section{$W$-entropy formula for geodesic flow on Wasserstein space}

Starting from Brenier's work  \cite{Br, BB} on the Monge-Kantorovich optimal transport problem with quadratic cost function,  Otto, Lott, McCann, Villani and Sturm \cite{Ot, OtV, LoV, Lo2, V1, V2, St1, St2, St3} among others have developed  the optimal transport theory. In particular, they developed an infinite dimensional Riemannian geometry and the theory of the gradient flow on the Wasserstein space over Euclidean space, compact Riemannian manifolds and metric measure spaces. The displacement convexity of the Boltzmann-Shannon entropy or the Renyi  entropy along geodesics on the Wasserstein space has been a key tool in \cite{LoV, Lo2, V1, V2, St1, St2, St3}  to introduce the notions of the upper bound of the dimension and the lower bound of the Ricci curvature on metric measure spaces. 
In \cite{MT}, McCann and Topping proved the contraction property of the $L^2$-Wasserstein distance between solutions of the backward heat equation
on closed manifolds equipped with the Ricci flow, which extends  previous results for the Fokker-Planck equation on Euclidean space  (due to Otto \cite{Ot})  and on complete Riemannian manifolds with suitable Bakry-Emery curvature condition (due to Sturm and von Renesse \cite{StR} ). See also \cite{T1, T2}. In \cite{Lo2}, Lott further proved two convexity results of the Boltzmann-Shannon 
type entropy along the geodesics on the Wasserstein space over closed manifolds equipped with the backward Ricci flow, 
which are  closely related to Perelman's result on the monotonicity of the $W$-entropy for the Ricci flow. In \cite{LL13b},  we extended Lott's convexity results to the Wasserstein space on compact Riemannian manifolds equipped with the backward Perelman Ricci flow.

Let $(M, g)$ be a complete Riemannian manifold equipped with a weighted volume measure $d\mu=e^{-\phi}dv$, where $\phi\in C^2(M)$ and $dv$ denotes the volume measure on $(M, g)$.  
The Boltzmann-Shannon entropy of the probability measure $\rho d\mu$ with respect to the reference measure $\mu$ is defined by
\begin{eqnarray*}
{\rm Ent}(\rho):= \int_M \rho \log \rho d\mu.
\end{eqnarray*}

Let $P_2(M, \mu)$ (resp. $P_2^\infty(M, \mu)$) be the Wasserstein 
space (reps. the smooth Wasserstein space) of all probability measures $\rho(x)d\mu(x)$ with density function (resp. with smooth density function) $\rho$ on $M$ such that $\int_M d^2(o, x)\rho(x)d\mu(x)<\infty$,  
where $d(o, \cdot)$ denotes the distance function from a fixed point $o\in M$. Following Otto \cite{Ot} and Lott \cite{Lo1, Lo2}, the tangent space $T_{\rho d\mu}P_2^\infty(M, \mu)$ is identified as follows 
\begin{eqnarray*}
T_{\rho d\mu}P_2^\infty(M, \mu)=\left\{s=\nabla_\mu^*(\rho \nabla f): f\in C^\infty(M), \ \ \int_M |\nabla f|^2\rho d\mu<\infty\right\},
\end{eqnarray*}
where
$\nabla_\mu^*$ denotes the $L^2$-adjoint of the Riemannian gradient $\nabla$ with respect to the weighted volume measure $d\mu$  on $(M, g)$.
For $s_i=\nabla_\mu^*(\rho\nabla f_i)\in T_{\rho d\mu} P_2^\infty(M, \mu)$,  we introduce Otto's infinite dimensional Riemannian metric on $P_2^\infty(M, \mu)$ as follows 
\begin{eqnarray*}
\langle \langle s_1, s_2\rangle\rangle:=\int_M \nabla f_1\cdot \nabla f_2 \rho d\mu,
\end{eqnarray*}
provided that 
\begin{eqnarray*}
\|s_i\|^2:=\int_M |\nabla f_i|^2\rho d\mu<\infty, \ \ \ i=1, 2.
\end{eqnarray*}
Let $T_{\rho d\mu}P_2(M, \mu)$ be the completion of $T_{\rho d\mu}P_2^\infty (M, \mu)$  equipped with Otto's  infinite dimensional  Riemannian metric. Then $P_2(M, \mu)$ is an infinite dimensional Riemannian manifold. 

By  Benamou and Brenier \cite{BB}, for any given $\mu_i=\rho_i d\mu\in P_2(M, \mu)$, $i=0, 1$, the $L^2$-Wasserstein distance between $\mu_0$ and $\mu_1$ coincides with the geodesic distance between $\mu_0$ and $\mu_1$ in $P_2(M, \mu)$ equippped with Otto's infinite dimensional Riemannian metric, i.e.,  
\begin{eqnarray*}
W_2^2(\mu_0, \mu_1)=\inf\limits\left\{{1\over 2}\int_0^1 |\nabla f(x, t)|^2\rho(x, t)d\mu(x): \partial_t \rho=\nabla_\mu^*(\rho \nabla f), \ \rho(0)=\rho_0, \ \rho(1)=\rho_1\right\}.
\end{eqnarray*}
By  \cite{Mc}, given $\mu_0=\rho(\cdot, 0)\mu, ~ \mu_1=\rho(\cdot, 1)\mu\in P_2^\infty(M, \mu)$ 
, it is known that there is a unique minimizing Wasserstein geodesic $\{\mu(t), t\in [0, 1]\}$  of the form $\mu(t) =(F_t)_*\mu_0$  
joining $\mu_0$ and $\mu_1$ in $P_2(M, \mu)$, where $F_t \in {\rm Diff}(M)$  
is given by $F_t(x) = \exp_x(-t \nabla f(\cdot, 0))$  for an appropriate Lipschitz function $f(\cdot, t)$. See also \cite{Lo1, Lo2}. 
 If the Wasserstein 
geodesic in $P_2(M, \mu)$  belongs entirely to $P_2^\infty(M, \mu)$,  then the geodesic flow $(\rho, f)\in T^*P_2^\infty(M, \mu)$ satisfies the transport equation  and the Hamilton-Jacobi equation 
\begin{eqnarray}
{\partial_t} \rho-\nabla_\mu^*(\rho \nabla f)&=&0,\label{TA}\\
{\partial_t}f+{1\over 2}|\nabla f|^2&=&0, \label{HJ}
\end{eqnarray}
with the boundary condition $\rho(0)=\rho_0$ and $\rho(1)=\rho_1$.   When  $\rho_0, f_0\in C^\infty(M)$,  defining $f(\cdot, t)\in C^\infty(M)$ by the Hopf-Lax solution
 \begin{eqnarray*}
f(x, t)=\inf\limits_{y\in M}\left(f_0(y)+{d^2(x, y)\over 2t}\right),\label{HLS}
 \end{eqnarray*} 
 and solving the transport equation $(\ref{TA})$ by the characteristic method,  it is known that $(\rho, f)$ satisfies $(\ref{TA})$ and $(\ref{HJ})$ with $\rho(0)=\rho_0$ and $f(0)=f_0$. See  \cite{V1}  Sect. 5.4.7. See also \cite{Lo1, Lo2}.
   In view of this, the transport equation $(\ref{TA})$ 
 and the Hamilton-Jacobi equation $(\ref{HJ})$ describe the geodesic flow on the cotangent bundle $T^*P_2^\infty(M, \mu)$ over the Wasserstein space $P_2(M, \mu)$. 
 Note that the Hamilton-Jacobi equation $(\ref{HJ})$ is also called the eikonal equation in geometric optics. 

The main result of this section is the following $W$-entropy formula for the geodesic flow on the Wasserstein space $P_2^\infty(M, \mu)$.

\begin{theorem}\label{MT2} (\cite{LL16, SLi15}) 
Let $(M, g)$  be a compact Riemannian manifold,  $\phi\in C^2(M)$, $d\mu=e^{-\phi}dv$.
Let $\rho: M\times [0, T]\rightarrow\mathbb{R}^+ $ and $f: M\times [0,T]\rightarrow \mathbb{R}$ be smooth solutions to the transport equation $(\ref{TA})$ 
and the Hamilton-Jacobi equation $(\ref{HJ})$. 
For any $m\geq n$, define the $H_m$-entropy and $W_m$-entropy for the geodesic flow $(\rho, f)$ on $T^*P^\infty_2(M, \mu)$ as follows
\begin{eqnarray*}
H_m(\rho, t)=-{\rm Ent}(\rho(t))-{m\over 2}\left(1+\log(4\pi t^2)\right),
\end{eqnarray*}
and
\begin{eqnarray*}
W_m(\rho, t)={d\over dt}(tH_m(\rho, t)).
\end{eqnarray*}
Then for all $t>0$, we have
\begin{eqnarray}
{d\over dt}W_m(\rho, t)&=&-t\int_M \left[\left|\nabla^2 f-{g\over t}\right|^2+Ric_{m,
n}(L)(\nabla f, \nabla f) \right]\rho d\mu\nonumber\\
& &\ \ \ \ \ \ \ \ \ \  -{t \over m-n}\int_M \left|\nabla \phi\cdot
\nabla f+{m-n\over t}\right|^2 \rho d\mu.\label{Wgeo}
\end{eqnarray}
In particular, if $Ric_{m, n}(L)\geq 0$, then $W_m(\rho, t)$ is decreasing in time $t$ along the geodesic flow on $T^*P^\infty_2(M, \mu)$.  
\end{theorem}

As a corollary of Theorem \ref{MT2}, we can recapture the following beautiful result due to Lott and Villani \cite{LoV, Lo2}.

\begin{corollary}\label{Th-LV}  (\cite{LoV, Lo2}) Let $(M, g, \phi)$ be a compact Riemannian manifold with $Ric_{m, n}(L)\geq 0$. Then $t{\rm Ent}(\rho(t))+mt\log t$ is convex in time $t$ along the geodesic on $P_2(M, \mu)$.
\end{corollary}

\section{Comparison between Theorem \ref{Th-A} and Theorem \ref{MT2}}

In this section,  we compare the $W$-entropy formula $(\ref{W-1})$ in Theorem \ref{Th-A} and the $W$-entropy formula $(\ref{Wgeo})$ in Theorem \ref{MT2}. 

\begin{itemize}

\item The $W$-entropy formula $(\ref{W-1})$ for the heat equation of the Witten Laplacian in Theorem \ref{Th-A} and the $W$-entropy formula $(\ref{Wgeo})$ for the geodesic flow on the Wasserstein space in Theorem \ref{MT2} have  similar expressions. Moreover, similarly to Corollary \ref{Th-LV}, from Theorem \ref{Th-A} we can derive the following
 
\begin{corollary} \label{Th-Li}  Let $(M, g, \phi)$ be a compact Riemannian manifold with $Ric_{m, n}(L)\geq 0$.  Then $t{\rm Ent}(u(t))+{m\over 2}t\log t$  is convex in time $t$ along the heat equation $\partial_t u=Lu$ on $M$. 
\end{corollary} 

\item By \cite{Li12, Li13}, Theorem \ref{Th-A} and  a rigidity theorem hold on complete Riemannian manifolds with bounded geometric condition and with the $CD(0, m)$-condition: $W_m(u, t)$ achieves its minimum at some $t=t_0>0$ if and only if $M=\mathbb{R}^n$, $m=n$, and $u(x, t)=u_m(x, t)={1\over (4\pi t)^{m\over 2}}e^{-{\|x\|^2\over 4t}}$ is the heat kernel of the heat equation $\partial_t u=\Delta u$  on $\mathbb{R}^m$. Note that, the Boltzmann-Shannon entropy of the Gussian heat kernel measure $u_m(x, t)dx$ 
is given by
$${\rm Ent}(u_m(t))=-{m\over 2}(1+\log(4\pi t)).$$
Thus the $H_m$-entropy for the heat equation of the Witten Laplacian is given by\footnote{Following Villani \cite{V1, V2}, we call $H_m(u(t))$ the {\it relative entropy} even though it is slightly different from the classical definition of the relative entropy in probability theory.}
\begin{eqnarray*}
H_m(u(t))={\rm Ent}(u_m(t))-{\rm Ent}(u(t)),
\end{eqnarray*}
and the $W_m$-entropy for the heat equation of the Witten Laplacian is given by the Boltzmann entropy formula
\begin{eqnarray}
W_m(u, t):={d\over dt}\left(t[{\rm Ent}(u_m(t))-{\rm Ent}(u(t))]\right).\label{Wm}
\end{eqnarray}
This gives a natural probabilistic interpretation of the $W$-entropy for the heat equation of the Witten Laplacian on Riemannian manifolds. 
See also Section $6$ for the probabilistic interpretation of the Perelman $W$-entropy for the Ricci flow.

\item On the other hand, when $m\in \mathbb{N}$, we can check that the following $(\rho_m, f_m)$ 
\begin{eqnarray*}
\rho_m(x, t)&=&{1\over (4\pi t^2)^{m/2}}e^{-{\|x\|^2\over 4t^2}},\\
f_m(x, t)&=&{\|x\|^2\over 2t},
\end{eqnarray*}
where $t>0, x\in \mathbb{R}^m$, is a  solution to the transport equation $(\ref{TA})$ and 
the Hamilton-Jacobi equation $(\ref{HJ})$ on $\mathbb{R}^m$ equipped with the standard Lebesgue measure, i.e.,  
\begin{eqnarray*}
{\partial_t} \rho+\nabla\cdot(\rho \nabla f)&=&0,\label{TAm}\\
{\partial_t}f+{1\over 2}|\nabla f|^2&=&0, \label{HJm}
\end{eqnarray*}
respectively. Moreover, the Boltzmann-Shannon entropy of the probability measure $\rho_m(t, x)dx$ (which equals to $u_m(t^2, x)dx$) is given by 
\begin{eqnarray*}
{\rm Ent}(\rho_m(t))=-{m\over 2}(1+\log(4\pi t^2)).
\end{eqnarray*}
Thus  we can reformulate  the $H_m$-entropy and the $W_m$-entropy for the geodesic flow on the Wasserstein space $P_2(M, \mu)$ as follows
\begin{eqnarray}
H_m(\rho(t))={\rm Ent}(\rho_m(t))-{\rm Ent}(\rho(t)), \label{NHW-2a}
\end{eqnarray}
and
\begin{eqnarray}
W_m(\rho, t):={d\over dt}\left(t[{\rm Ent}(\rho_m(t))-{\rm Ent}(\rho(t))]\right).\label{NHW-3a}
\end{eqnarray}

\item The relative entropy  $H_m(\rho(t))$ defined by $(\ref{NHW-2a})$ is the difference between the Boltzmann-Shannon entropy of the probability measure $\rho(t) d\mu$ on $(M,\mu)$ and the Boltzmann-Shannon entropy of the reference model $\rho_m(t)dx$ on $(\mathbb{R}^m, dx)$, and $W_m(\rho, t)$ defined by $(\ref{NHW-3a})$  is the time derivative of $tH_m(\rho(t))$. In \cite{LL16, SLi15}, similarly to the case of Theorem \ref{Th-A},  we extended  the $W$-entropy formula $(\ref{Wgeo})$ in Theorem \ref{MT2}  to complete Riemannian manifolds with bounded geometry condition. In view of this,  we proved  that the rigidity model for the $W$-entropy for the geodesic flow on the Wasserstein space $P_2^\infty(M, \mu)$ over 
complete Riemannian manifolds with the $CD(0, m)$-condition is  $M=\mathbb{R}^n$, $m=n$, $\rho=\rho_m$ and $f=f_m$.

\end{itemize}

\section{Langevin deformation of geometric flows on Wasserstein space}

We can raise a natural question how to understand  the similarity between the $W$-entropy formulas in Theorem \ref{Th-A} and Theorem \ref{MT2}. Can we pass through one of them to another one? One possible approach to answer this question is to use the vanishing viscosity limit method from the heat equation to the Hamilton-Jacobi equation. However,  it seems that one cannot  easily use this approach to pass through the $W$-entropy formula for the heat equation of the Witten Laplacian to the $W$-entropy formula for the geodesic flow on the Wasserstein space. 

Inspired by J.-M.Bismut's works (see \cite{Bis05, Bis10}) on the deformation of hypoelliptic Laplacians on the cotangent bundle over Riemannian manifolds,  which interpolates the usual Laplacian on the underlying Riemannian manifold $M$ and the Hamiltonian vector field which generates the geodesic flow on the cotangent bundle over $M$,  we introduced in  \cite{LL16, SLi15}  the Langevin 
deformation of geometric flows on the cotangent bundle of the Wasserstein space over compact Riemannian manifolds

More precisely, for $c\in (0, \infty)$,  let $(\rho, f)$ be smooth solution to the following equations

\begin{eqnarray}
\partial_t \rho-\nabla_\mu^*(\rho\nabla f)&=&0,\label{flow1}\\
c^2\left(\partial_t f+{1\over 2}|\nabla f|^2\right)&=&-f+V'(\rho). \label{flow2}
\end{eqnarray}
where $V\in C^\infty ((0, \infty),  \mathbb{R})$. Eq.~$(\ref{flow1})$ is indeed the transport equation $(\ref{TA})$, while Eq.~$(\ref{flow2})$ can be viewed as the 
the Langevin equation on $T^*P_2(M, \mu)$. 
When $c\rightarrow \infty$,  Eq.~$(\ref{flow2})$ implies that $f$ should satisfies the Hamilton-Jacobi equation $(\ref{HJ})$. In this case, $(\rho, f)$ is indeed the geodesic flow on $T^*P_2(M, \mu)$. 
On the other hand, when $c=0$,   Eq.~$(\ref{flow2})$ implies that $f=V'(\rho)$. In this case,  $\rho$ is the backward gradient flow of $U(\rho)=\int_M V(\rho)d\mu$ on  
$P_2(M, \mu)$ equipped with Otto's infinite dimensional Riemannian metric
\begin{eqnarray}
\partial_t \rho=\nabla_\mu^*\left(\rho\nabla V'(\rho)\right). \label{gradflow}
\end{eqnarray}

In \cite{LL16, SLi15}, we proved the existence and uniqueness of the Langevin deformation between the backward gradient flow of the Boltzmann-Shannon entropy  ${\rm Ent}(\rho)=\int_M \rho \log \rho d\mu$ (respectively,  the Renyi entropy $U(\rho)=
{1\over m-1}\int_M \rho^{m}d\mu$ for $m>1$), which is 
the backward heat equation $\partial_t \rho=-L\rho$ (respectively, the backward porous medium equation $\partial_t \rho=-L\rho^m$) of the Witten Laplacian on $M$, and the geodesic flow on the  cotangent bundle of the smooth Wasserstein space $P_2^\infty(M, \mu)$ over compact Riemannian manifolds. Moreover, 
we  proved an extension of the $W$-entropy formula along the Langevin deformation of geometric flows. The rigidity models are also proposed for the Langevin deformation.  Due to the limit of the length of the paper, we refer  
the reader to \cite{LL16, SLi15} for the details of these results.  

\section{The $W$-entropy, statistical mechanics and probability theory}

In \cite{P1}, Perelman gave a heuristic
interpretation for the $W$-entropy using statistical mechanics. Recall that the partition function for the canonical ensemble at temperature $\beta^{-1}$  is given by $Z_\beta=\int_{\mathbb{R}} e^{-\beta E}d\omega(E)$, where $d\omega(E)$
denotes the ``density of states'' measure, whose physical meaning is
the number of microstates with energy levels in the range $[E,
E+dE]$. The average of the energy with respect to the Gibbs measure $dP(E)={e^{-\beta E} \over Z_\beta} d\omega(E)$  is 
$$\langle E\rangle=-{\partial\over \partial \beta}\log Z_\beta,$$
and the Boltzmann entropy $S$ satisfies the Boltzmann entropy formula 
\begin{eqnarray*}
S=\log Z_\beta-\beta {\partial \over \partial \beta}\log Z_\beta.
\end{eqnarray*}
Equivalently, letting $\tau=\beta^{-1}$, then
\begin{eqnarray}
S= {\partial \over \partial \tau}(\tau \log Z_\beta).\label{S1}
\end{eqnarray}
The fluctuation of the energy is given by 
$$
\sigma:=\langle (E-\langle E\rangle)^2\rangle ={\partial^2 \over
\partial \beta^2} \log Z_\beta,$$
and the derivative of the Boltzmann entropy with respect to $\beta$ satisfies
\begin{eqnarray*}
{\partial S\over \partial \beta}=-\beta \sigma.\label{S2}
\end{eqnarray*}

Let
$(M, g(\tau))$ be a family of closed Riemannian manifolds,
$dm(\tau)=(4\pi\tau)^{-n/2}e^{-f(\tau)}dv_{g(\tau)}$ a probability
measure on $(M, g(\tau))$, where $g(\tau)$ satisfies the backward Ricci flow equation $\partial
_\tau g=2Ric$,  $f(\tau)$ satisfies
the heat equation $\partial_\tau f=\Delta f-|\nabla f|^2+R-{n\over
2\tau}$ and $\tau=T-t$. Assume that  there is a canonical ensemble with a ``density of state
measure'' $d\omega(E)$ such that the partition function
$Z_\beta=\int_{\mathbb{R}} e^{-\beta E}d\omega(E)$ is given by
\begin{eqnarray}
\log Z_\beta=\int_M \left(-f+{n\over 2}\right)dm,\label{logZ}
\end{eqnarray}
where $\beta={1\over \tau}$, and the backward time $\tau=T-t$ is regarded as the
temperature. Then, using the above formulas in statistical mechanics,  Perelman \cite{P1} formally derived that 
\begin{eqnarray*}
\langle E\rangle&=&=-\tau^2\int_M \left(R+|\nabla f|^2-{n\over
2\tau}\right)dm,\\
S&=&-\int_M \left(\tau(R+|\nabla f|^2)+f-n\right)dm,\\
\sigma&=&2\tau^4\int_M\left|Ric+\nabla^2f-{g\over
2\tau}\right|^2dm.
\end{eqnarray*}
This yields 
\begin{eqnarray*}
W(g, f, \tau)=-S,\label{W1}
\end{eqnarray*}
and
\begin{eqnarray*}
{d\over dt}W(g, f, \tau)=2 \int_M \tau\left|Ric+\nabla^2 f-{g\over
2\tau}\right|^2dm.\label{W2}
\end{eqnarray*}
This gives an interpretation of the $W$-entropy $(\ref{entropy-1})$ for the backward Ricci flow by Boltzmann's entropy formula $(\ref{S1})$. However,  the problem whether there is a canonical
ensemble  with a ``density of states'' measure $d\omega(E)$ such that
the partition function $Z_\beta=\int_{\mathbb{R}} e^{-\beta E}d\omega(E)$
satisfies Perelman's requirement $(\ref{logZ})$  remains open.  See \cite{Li12b} for further discussion on this issue.

In \cite{Li12, Li12b}, the second author of this paper gave a probabilistic interpretation of
Perelman's $W$-entropy for the Ricci flow. Observing that
\begin{eqnarray*}
\log Z_\beta={\rm Ent}(u(\tau))+{n\over 2}(1+\log (4\pi \tau)),
\end{eqnarray*}
where
\begin{eqnarray*}
{\rm Ent}(u(\tau))=\int_M u\log udv=-\int_M \left(f+{n\over 2}\log(4\pi
\tau)\right){e^{-f}\over (4\pi\tau)^{n/2}}dv
\end{eqnarray*}
is the Boltzmann-Shannon entropy of the heat kernel measure
$dm=u(\tau)dv_{g(\tau)}$ with respect to the volume measure $dv_{g(\tau)}$ on $(M, g(\tau))$, where $u={e^{-f}\over (4\pi\tau)^{n/2}}$. On the other
hand, let
$$u_n(x, \tau)={e^{-{\|x\|^2\over 4\tau}}\over (4\pi \tau)^{n/2}}, \ \ \ \forall x\in \mathbb{R}^n, \tau>0$$
be the Gaussian heat kernel on $\mathbb{R}^n$. Then it is well-known
that the Boltzmann-Shannon entropy of the Gaussian measure $d\gamma_n(\tau, x)=u_n(\tau, x)dx$ with respect to the Lebesgue measure is given by
\begin{eqnarray*}
{\rm Ent}(u_n)=-{n\over 2}(1+\log (4\pi \tau)).
\end{eqnarray*}
Hence
$$\log Z_\beta={\rm Ent}(u(\tau))-{\rm Ent}(u_n(\tau))$$
 is the difference of  the Boltzmann-Shannon entropy of the heat
kernel measure $dm=u(\tau)dv_{g(\tau)}$ on $(M, g(\tau))$ and the
Boltzmann-Shannon entropy of the Gaussian measure $\gamma_n$ on
$\mathbb{R}^n$. In view
of this, we have the following probabilistic interpretation of the
$W$-entropy for the Ricci flow
\begin{eqnarray}
W(g, f, \tau):={d\over d\tau}(\tau [{\rm Ent}(u_n(\tau))-{\rm Ent}(u(\tau))]). \label{WEnt}
\end{eqnarray}

Similarly to $(\ref{WEnt})$, we can give the probabilistic interpretation of the $W$-entropy for the heat equation of the Witten Laplacian on complete Riemannian manifolds with 
the $CD(0, m)$-condition. See  $(\ref{Wm})$ in Section $5$. By $(\ref{NHW-2a})$ and $(\ref{NHW-3a})$, we have the similar probabilistic interpretation of the $W$-entropy for 
the geodesic flow on the Wasserstein space over Riemannian manifolds with 
the $CD(0, m)$-condition. See \cite{Li12}. 

It is natural and interesting to ask the question whether we can give a probabilistic interpretation of the $W$-entropy for the heat equation of the Witten Laplacian on complete Riemannian manifolds with 
the $CD(K, m)$ and $CD(K, \infty)$-conditions. 
This question is closely related to the question whether there exist the rigidity models of the $W$-entropy for the heat equation of the Witten Laplacian on complete Riemannian manifolds with the $CD(K, m)$ and $CD(K, \infty)$-conditions. 

In our recent paper \cite{LL17a}, we gave a probabilistic interpretation of the $W_{m, K}$-entropy for the heat equation of the Witten Laplacian on complete Riemannian manifolds with the $CD(K, m)$-condition. More precisely,  let $m\in \mathbb{N}$, $M=\mathbb{R}^m$, $g_0$ the Euclidean metric, $\phi_K(x)=-{K\|x\|^2\over 2}$ and $d\mu_K(x)=e^{K\|x\|^2\over 2}dx$, where $K\in \mathbb{R}$.  Then $\nabla\phi_K(x)=-Kx$, and $\nabla^2\phi_K=-K{\rm Id}_{\mathbb{R}^m}$. We consider the Ornstein-Ulenbeck operator on $\mathbb{R}^m$ given by
$$L=\Delta+Kx\cdot \nabla.$$
Note that  $(\mathbb{R}^m, g_0, \phi_{K})$ is a complete shrinking Ricci soliton, i.e., $Ric(L)=-Kg_0$. The  Ornstein-Ulenbeck diffusion process on $\mathbb{R}^m$ is the solution to the Langevin SDE
$$dX_t=\sqrt{2} dW_t+KX_tdt, \ \ \ X_0=x.$$
It is well-known that  the law of $X_t$ is Gaussian  $N\left(e^{Kt}x, {e^{2Kt}-1\over K}{\rm Id}\right)$, and the heat kernel of $X_t$ with respect to the Lebesgue measure on $\mathbb{R}^m$  is given by
\begin{eqnarray*}
u_{m, K}(x, y, t)=\left({K\over 2\pi (e^{2Kt}-1)}\right)^{m/2}\exp \left(-{K|y-e^{Kt}x|^2\over 2(e^{2Kt}-1)}\right).
\end{eqnarray*}
By direct calculation, the relative Boltzmann-Shannon entropy of the law of $X_t$ with respect to the Lebesgue measure on $\mathbb{R}^m$  is given by
\begin{eqnarray*}
{\rm Ent}(u_{m, K}(x, y, t)| dy)=-{m\over 2}\left(1+\log(4\pi \sigma_K^2(t))\right),
\end{eqnarray*}
where $\sigma_K^2(t)={e^{2Kt}-1\over 2K}$. When $t\rightarrow 0$, we have
\begin{eqnarray*}
{\rm Ent}(u_{m, K}(x, y, t) | dy)=-{m\over 2}\left(1+\log (4\pi t)+Kt+{K^2t^2\over 6}\right)+O(t^4).
\end{eqnarray*}

Thus, when $t\rightarrow 0^+$, the second term in the definition formula $(\ref{HmK})$ of the $H_{m, K}$-entropy is  asymptotically equivalent (at the order $O(t^4)$)  to the Boltzmann-Shannon entropy of the heat kernel at time $t$ of the Ornstein-Uhlenbeck operator on $\mathbb{R}^m$ with respect to the Lebesgue measure on $\mathbb{R}^m$.  That is to say,  when $t\rightarrow 0^+$, we have
\begin{eqnarray*}
H_{m, K}(u(t))={\rm Ent}(u_{m, K}(t)|dy)-{\rm Ent}(u(t)|\mu)+O(t^4),
\end{eqnarray*}
and
\begin{eqnarray*}
W_{m, K}(u(t))={d\over dt}\left(tH_{m, K}(u(t))\right).
\end{eqnarray*}

To end this paper, let us mention that, in a forthcoming paper \cite{KL16}, Kuwada and the second author of this paper prove an analogue of the $W$-entropy monotonicity theorem on metric measure spaces with the so-called $RCD(0, N)$-condition.

\medskip

\noindent{\bf Acknowledgement}.  We would like to thank Professors S. Aida, 
D. Bakry, J.-M. Bismut, D. Elworthy, K. Kuwae, M. Ledoux, N. Mok, K.-T. Sturm, A. Thalmaier, F.-Y. Wang  and Dr. Yuzhao Wang for their interests and 
helpful discussions 
during various stages of this work.

\begin{flushleft}
\medskip\noindent

Songzi Li, School of Mathematical Science, Beijing Normal University, No. 19, Xin Jie Kou Wai Da Jie, 100875, China, Email: songzi.li@bnu.edu.cn

\medskip

Xiang-Dong Li, Academy of Mathematics and Systems Science, Chinese
Academy of Sciences, 55, Zhongguancun East Road, Beijing, 100190, China, 
E-mail: xdli@amt.ac.cn
\\
and
\\
School of Mathematical Sciences, University of Chinese Academy of Sciences, Beijing, 100049, China
\end{flushleft}

\end{document}